\documentclass{amsart}
\usepackage{graphicx}
\usepackage{amsmath}
\usepackage{amssymb}
\usepackage[all]{xy}
\usepackage{amscd}
\usepackage{bbold}

\newcommand{\Id}{\operatorname{Id}}
\newcommand{\Ktac}{\operatorname{K_{tac}}}
\newcommand{\Hom}{\operatorname{Hom}}
\renewcommand{\H}{\operatorname{H}}
\renewcommand{\Im}{\operatorname{Im}\nolimits}
\newcommand{\Ker}{\operatorname{Ker}\nolimits}
\newcommand{\Ann}{\operatorname{Ann}\nolimits}
\newcommand{\Ext}{\operatorname{Ext}\nolimits}

\newcommand{\Cone}{\operatorname{Cone}\nolimits}

\newcommand{\proj}{\operatorname{proj}}
\newcommand{\Supp}{\operatorname{Supp}\nolimits}
\newcommand{\Rank}{\operatorname{rank}\nolimits}

\newtheorem*{prop1.1}{Proposition 1.1}
\newtheorem*{prop2.1}{Proposition 2.1}
\newtheorem*{prop2.2}{Proposition 2.2}
\newtheorem*{prop2.3}{Corollary 2.3}
\newtheorem*{prop2.4}{Corollary 2.4}
\newtheorem*{prop2.5}{Corollary 2.5}
\newtheorem*{prop2.6}{Proposition 2.6}
\newtheorem*{lem4.1}{Lemma 4.1}
\newtheorem*{lem4.2}{Lemma 4.2}
\newtheorem*{lem4.3}{Lemma 4.3}
\newtheorem*{thm3.3}{Theorem 3.3}
\newtheorem*{thm3.5}{Theorem 3.5}
\newtheorem*{thm4.4}{Theorem 4.4}
\newtheorem*{thm5.1}{Theorem 5.1}
\newtheorem*{thm6.1}{Theorem 6.1}
\newtheorem*{thm6.2}{Theorem 6.2}
\newtheorem*{thm6.3}{Theorem 6.3}
\newtheorem*{thm6.4}{Theorem 6.4} 
\newtheorem*{thm6.5}{Theorem 6.5}

\theoremstyle{definition}
\newtheorem*{def3.1}{Definition 3.1}
\newtheorem*{def3.4}{Definition 3.4}
\newtheorem*{def3.2}{Definition 3.2}
\newtheorem*{example1}{Example 7.1}
\newtheorem*{example2}{Example 7.2}
\newtheorem*{Remark}{Remark}
\begin{document}
	
\title{Support and Rank Varieties of Totally Acyclic Complexes}

\author{Nathan T. Steele}

\date{\today}



\begin{abstract} 
Support and rank varieties of modules over a group algebra of an elementary abelian $p$-group have been well studied. In particular, Avrunin and Scott showed that in this setting, the rank and support varieties are equivalent. Avramov and Buchweitz proved an analogous result for pairs of modules over arbitrary commutative local complete intersection rings. In this paper we study support and rank varieties in the triangulated category of totally acyclic chain complexes over a complete intersection ring and show that these varieties are also equivalent.  
\end{abstract}

\maketitle

\section*{Introduction}
The main goal of this paper is to explore a theory of support and rank varieties in the triangulated category of totally acyclic chain complexes $\Ktac(R)$ over a complete intersection ring $R=Q/(f_1,\dots,f_c)$. In \cite[1.1]{Avrunin}, Avrunin and Scott showed that for a module $M$ over a group algebra of an elementary abelian $p$-group, its support variety $V(M)$ is equivalent to its rank variety $W(M)$. Avramov and Buchweitz \cite[2.5]{Avramov:2000} proved a similar result for pairs of modules $M,N$ over a complete intersection ring and they simply refer to both sets as the support variety $V(M,N)$ . We give similar definitions of rank and support varieties in the category of totally acyclic chain complexes and prove that these varieties are equal.  

The proof of the equivalence of support varieties of modules over a complete intersection ring requires a theorem that states for any pair of $R$-modules $M,N$, the graded $R$-module $\Ext_Q^*(M,N)$ is finitely generated if and only if the graded $R[\chi_1,\dots,\chi_c]$-module $\Ext_R^*(M,N)$ is finitely generated, where $R[\chi_1,\dots,\chi_c]$ is a ring of cohomology operators. The `if' direction of this theorem was proved by Avramov, Gasharov, and Peeva \cite[4.2]{Avramov:1997} and the `only if' direction was proved by Gulliksen \cite[3.1]{Gulliksen:1974} cf. \cite[2.1]{Avramov:1989}. We take advantage of the triangulated structure of $\Ktac(R)$ and of a pair of adjoint functors $T:\Ktac(R)\rightarrow \Ktac(Q)$ and $S:\Ktac(Q)\rightarrow \Ktac(R)$ to prove a similar finite generation result for pairs of totally acyclic complexes. We then use a slightly modified version of the proof by Avramov and Buchweitz \cite[2.5]{Avramov:2000} to prove the equivalence of varieties.

\section{Totally Acyclic Chain Complexes}

Let $R$ be an associative ring, $C$ and $D$ be chain complexes of $R$-modules, and $f,f':C\rightarrow D$ be chain maps. We say $f$ is \emph{homotopic} to $f'$ if there exist maps $\lambda_n:C_n\rightarrow D_{n+1}$ such that $\psi_n - \psi'_n = d^D_{n+1}\lambda_n + \lambda_{n-1}d^C_n$ and we write $\psi \sim \psi'$.
$$ 
\xymatrix{
C:\cdots \ar[rr]^{d_{n+2}^C}  && C_{n+1} \ar[rr]^{d_{n+1}^C} \ar[dd]_{\psi_{n+1}}^{\psi'_{n+1}} && C_n \ar[rr]^{d_n^C} \ar[dd]_{\psi_{n}}^{\psi'_{n}} \ar@{.>}[ddll]_{\lambda_{n}} && C_{n-1}\ar[rr]^{d_{n-1}^C} \ar[dd]_{\psi_{n-1}}^{\psi'_{n-1}} \ar@{.>}[ddll]_{\lambda_{n-1}} && \cdots \\
\\
D:\cdots \ar[rr]^{d_{n+2}^D} && D_{n+1} \ar[rr]^{d_{n+1}^D} && D_n \ar[rr]^{d_n^D} && D_{n-1}\ar[rr]^{d_{n-1}^D} && \cdots
} $$
 If there exist maps $\psi:C\rightarrow D$ and $\zeta:D\rightarrow C$ such that $\psi\circ \zeta \sim \Id_D$ and $\zeta \circ \psi \sim \Id_C$, then we say $C$ and $D$ are homotopically equivalent and we write $C\simeq D$.  

We denote by $\Hom_R(C,D)$ the set of homotopy equivalence classes of degree zero chain maps from $C$ to $D$. The complex $C$ is called \emph{totally acyclic} if each $C_i$ is a finitely generated free $R$-module, $C$ is exact, and $\Hom_R(C,R)$ is exact, i.e. 
$$\H(C)=0=\H(\Hom_R(C,R)).$$

The set of totally acyclic chain complexes forms a triangulated subcategory of the homotopy category of chain complexes. We denote this subcategory by $\Ktac(R)$ and it has distinguished triangles of the form
$$
\begin{CD}
C @>t>> D @>>> \Cone(t) @>>> \Sigma C 
\end{CD}
$$
where $t$ is a zero degree chain map from $C$ to $D$. The shift functor $\Sigma$ simply moves every module in the complex to the left, i.e $(\Sigma^i C)_n = C_{n-i}$. For simplicity we will write $\Sigma C=\Sigma^1 C$. The differentials of the shifted complex are given by $d_n^{\Sigma^i C} =  (-1)^i(\sigma_{C}^i)_{n-i-1} d_{n-i}^C (\sigma^{i}_{C})_{n}^{-1}$ where $\sigma^i_{C}$ is the natural isomorphism $\sigma^i_{C}: C \rightarrow \Sigma^i C$. Notice that $\sigma_{\Sigma^iC}^{-i} = (\sigma_C^i)^{-1}$. 

The \emph{Mapping Cone of t} is the chain complex $\Cone(t) = \Sigma C\oplus D$ with differentials given by
	$$d^{\Cone(t)}_n =  \left(\begin{array}{cc} d^{\Sigma C}_n & 0 \\  t_{n-1} (\sigma^{1}_{C})_{n-1}^{-1}  & d^D_{n}  \end{array}\right).$$

\subsection*{Complete Resolutions}

 Let $M$ be an $R$-module and let $P$ be a projective resolution of $M$. A \emph{Complete Resolution of $M$} is any totally acyclic chain complex $C$ such that the truncated complex $C_{>i}:\cdots \rightarrow C_{i+2}\rightarrow C_{i+1} \rightarrow 0$ is homotopically equivalent to $P_{>i}:\cdots \rightarrow P_{i+2} \rightarrow P_{i+1} \rightarrow 0$ for some $i$.

If we let $P'$ be a projective resolution of a module $N$ and let $\bar{\rho_0}:M\rightarrow N$ and $\bar{\mu_0}:N\rightarrow M$ be homomorphisms, then the comparison theorem \cite[6.16]{Rotman} implies that there exist chain maps $\rho:P \rightarrow P'$ and $\mu:P'\rightarrow P$. Furthermore, any chain map from $P$ to $P'$ lifting $\bar{\rho_0}$ is homotopic to $\rho$ and any chain map from $P'$ to $P$ lifting $\bar{\mu_0}$ is homotopic to $\mu$.

\begin{prop1.1} 
 \emph{(Avramov-Martsinkovsky \cite[5.3]{Avramov:2002})} Complete resolutions are unique up to homotopy equivalence. 
 \end{prop1.1}
The following is an alternate proof to that given by Avramov and Martsinkovsky.
 \begin{proof}
If $C,C'$ are complete resolutions of $M,N$ respectively, then $C_{>i}\simeq P_{>i}$ and $C'_{>k}\simeq P'_{>k}$ for some $i,k$. Thus there exist maps $h_{>i}:C_{>i}\rightarrow P_{>i}$, $h_{>k}':C'_{>k}\rightarrow P'_{>k}$, $j_{>i}:P_{>i}\rightarrow C_{>i}$, and $j_{>k}':P'_{>k}\rightarrow C'_{>k}$ such that $h_{>i}\circ j_{>i}\sim\Id_{P_{>i}}$, $j_{>i} \circ h_{>i}\sim \Id_{C_{>i}}$, $h_{>k}'\circ j_{>k}'\sim\Id_{P_{>k}'}$, and $j_{>k}' \circ h_{>k}'\sim \Id_{C_{>k}'}$.  

We can set $l=\max(i,k)$ so that $C_{>l}\simeq P_{>l}$ and $C'_{>l}\simeq P'_{>l}$. Hence for any homomorphisms $\bar{\rho_0}:M\rightarrow N$ and $\bar{\mu_0}:N\rightarrow M$ with $\rho:P \rightarrow P'$ and $\mu:P'\rightarrow P$ as above, there exist chain maps $\gamma_{>l}:C_{>l} \rightarrow C'_{>l}$ and $\delta_{>l}:C'_{>l}\rightarrow C_{>l}$ where $\gamma_n = j'_n \circ \rho_n \circ h_n $ and $\delta_n = j_n \circ \mu_n \circ h'_n$ for $n\geq l$.
Now consider the following sequence. $$\cdots \rightarrow \Hom(C_l,C_l')\rightarrow \Hom(C_{l+1},C_l')\rightarrow \Hom(C_{l+2},C_l') \rightarrow \cdots $$
This sequence is exact since $C$ is totally acyclic. Now consider following the diagram.
$$\xymatrix{ 
\cdots \ar[r]^{d_{l+3}^C}  & C_{l+2} \ar[r]^{d_{l+2}^C} \ar[d]_{\gamma_{l+2}} & C_{l+1} \ar[r]^{d_l^C} \ar[d]_{\gamma_{l+1}}  & C_{l}\ar[r]^{d_{l}^C} \ar@{.>}[d]_{\gamma_{l}}  & \cdots \\
\cdots \ar[r]^{d_{l+3}^{C'}} & C'_{l+2} \ar[r]^{d_{l+2}^{C'}} & C'_{l+1} \ar[r]^{d_{l+1}^{C'}} & C'_{l}\ar[r]^{d_{l}^{C'}} & \cdots}
$$
We would like to show that $\gamma_n$ exists for $n\leq l$ such that the diagram commutes. We know that the diagram commutes for $C_{>l}$, thus 
$$d_{l+1}^{C'}\circ \gamma_{l+1} \circ d_{l+2}^C = d_{l+1}^{C'} \circ d_{l+2}^{C'} \circ \gamma_{l+2} = 0$$ 
Therefore we have 
$$d^{C'}_{l+1}\circ \gamma_{l+1} \in \ker (\Hom(d^C_{l+2},C_l')) = \Im (\Hom(d^C_{l+1},C_l')).$$ 
Then there exists $\gamma_l\in \Hom(C_l,C_l')$ such that $d^{C'}_{l+1}\circ \gamma_{l+1} = \gamma_l \circ d^C_{l+1}$. We can now apply induction to construct $\gamma_n:C_n \rightarrow C'_n$ for $n\leq l$ such that $d^{C'}_{n+1}\circ \gamma_{n+1} = \gamma_n \circ d^C_{n+1}$  and similarly for $\delta_n:C'_n\rightarrow C_n$.

Thus we now have chain maps $\gamma:C\rightarrow C'$ and $\delta:C'\rightarrow C$. If $M=N$ and $\rho_0=\mu_0= \Id_M$, then $\rho \circ \mu \sim \Id_{p'}$ and $\mu \circ \rho \sim \Id_P$. Therefore $\gamma \circ \delta \sim \Id_C$ and $\delta \circ \gamma \sim \Id_{C'}$, hence $C\simeq C'$.
\end{proof}
\section{Eisenbud and Cohomology Operators in $\Ktac(R)$}
From this point forward we will let $R=Q/(f_1,\dots,f_c)$ where $(Q,m,\mathbb{k})$ is a local ring and $f=f_1,\dots,f_c$ is a regular sequence in $m^2$. We will also assume that the residue field $\mathbb{k}$ is algebraically closed. For pairs of modules over a complete intersection ring, there are multiple equivalent methods of constructing cohomology operators. To construct similar operators in $\Ktac(R)$, we will use Eisenbud's \cite{Eisenbud:1980} method of lifting chain complex differentials to create degree $2$ chain maps.   

Given an $R$-complex $C$, we define $\widetilde{C}$ to be a sequence of $Q$-free module homomorphisms such that
$$ \widetilde{C}\otimes_Q R \simeq C$$	
and we call $\widetilde{C}$ a \emph{lift} of $C$. Note that $\widetilde{C}$ is not necessarily a chain complex. In particular, 
$$\widetilde{d}_{n-1}^C \circ \widetilde{d}_n^C = \sum\limits_{k=1}^c f_k\widetilde{t}_{k,n}$$
where $\widetilde{t}_{k,n}:\widetilde{C}_n\rightarrow \widetilde{C}_{n-2}$. We will write $\widetilde{t}_k$ for the family of maps $\{\widetilde{t}_{k,n}\}_{n\in \mathbb{Z}}$.

\begin{prop2.1}
\emph{(Eisenbud \cite[1.1]{Eisenbud:1980})} The map $t_k = \sigma^2_C(\widetilde{t}_k \otimes_Q R):C\rightarrow \Sigma^2C$ is a morphism.
\end{prop2.1}
Now let $D$ be another chain complex with 
$$\widetilde{d}_{n-1}^D \circ \widetilde{d}_n^D = \sum\limits_{k=1}^c f_k\widetilde{s}_{k,n} \;\; \text{and} \;\; s_k = \sigma^2_D (\widetilde{s}_k \otimes_Q R).$$
\begin{prop2.2}
\emph{(Eisenbud \cite[1.3]{Eisenbud:1980})} 
If $\psi:C\rightarrow D$ is a chain map, then 
$$\psi \circ (\sigma_{C}^{2})^{-1} \circ t_k \sim (\sigma_{D}^{2})^{-1} \circ s_k \circ \psi.$$
\end{prop2.2}
We will now give a reorganized and more detailed version of Eisenbud's original proof.
\begin{proof}
Since $\psi$ is an $R$-complex chain map, lifting it to $Q$ yields
$$\widetilde{\psi}\widetilde{d}^C -\widetilde{d}^D\widetilde{\psi} = \sum\limits_{j=1}^c f_j \widetilde{\tau}_j.$$
We can also lift the maps $t_k$ and $s_k$ to $Q$ and use the facts that $(\widetilde{d}^{C})^2=\sum\limits_{j=1}^c f_j\widetilde{t}_j$ and $(\widetilde{d}^{D})^2=\sum\limits_{j=1}^c f_j\widetilde{s}_j$ to get
$$f_k(\widetilde{\psi}\widetilde{t}_k-\widetilde{s}_k\widetilde{\psi}) =  \widetilde{\psi} \left((\widetilde{d}^{C})^2-\sum\limits_{j \neq k} f_j\widetilde{t}_j\right)-\left((\widetilde{d}^{D})^2-\sum\limits_{j \neq k} f_j\widetilde{s}_j\right)\widetilde{\psi}.$$
We will now commute the $\widetilde{\psi}$ with one copy each of $\widetilde{d}^C$ and $\widetilde{d}^D$ and subtract like terms to get
$$f_k(\widetilde{\psi}\widetilde{t}_k-\widetilde{s}_k\widetilde{\psi}) = \sum\limits_{j \neq k} f_j \widetilde{\psi}(\widetilde{s}_j-\widetilde{t}_j)+\sum\limits_{j=1}^c f_j(\widetilde{d}^D\widetilde{\tau}_j +\widetilde{\tau}_j\widetilde{d}^C).$$
Taking $f_k(\widetilde{d}^D\widetilde{\tau}_k +\widetilde{\tau}_k\widetilde{d}^C)$ to the left yields
$$f_k\left(\widetilde{\psi}\widetilde{t_k}-\widetilde{s_k}\widetilde{\psi}-(\widetilde{d}^D\widetilde{\tau}_k +\widetilde{\tau}_k\widetilde{d}^C)\right) = \sum\limits_{j \neq k} f_j \widetilde{\psi}(\widetilde{s}_j-\widetilde{t}_j)+\sum\limits_{j \neq k} f_j(\widetilde{d}^D\widetilde{\tau}_j +\widetilde{\tau}_j\widetilde{d}^C).$$
Thus we see that $$\Im\left(\widetilde{\psi}\widetilde{t}_k-\widetilde{s}_k\widetilde{\psi}-(\widetilde{d}^D\widetilde{\tau}_k +\widetilde{\tau}_k\widetilde{d}^C)\right) \subseteq (f_1,\dots,f_c)\Sigma^iD$$ 
and so tensoring down to $R$ we have
$$\psi  (\sigma_{C}^{2})^{-1}  t_k - (\sigma_{D}^{2})^{-1}  s_k  \psi= d^D \tau_k +\tau_k d^C.$$
We can consider the map $\tau_k$ as a homotopy map so that $\psi  (\sigma_{C}^{2})^{-1}  t_k - (\sigma_{D}^{2})^{-1}  s_k  \psi \sim 0$ as desired. 
\end{proof}
The next two corollaries show that the $t_k$'s are well defined in $\Ktac(R)$. 

\begin{prop2.3}\emph{(Eisenbud \cite[1.4, 1.5]{Eisenbud:1980})} The $t_k$'s are uniquely determined up to homotopy and $t_kt_j \sim t_jt_k$. 
\end{prop2.3}

\begin{prop2.4}
If $C$ and $D$ are homotopically equivalent with $\psi:C\rightarrow D$ and $\zeta:D\rightarrow C$ such that $\psi \zeta \sim \Id_D$ and $\zeta \psi \sim \Id_C$, then $t_k \sim \zeta s_k \psi$ and $s_k \sim \psi t_k \zeta$.
\end{prop2.4}
\begin{proof}
By Proposition 2.2, $\psi t_k \sim s_k \psi$ and thus $t_k \sim \zeta \psi t_k \sim \zeta s_k \psi$. Similarly $s_k \sim \psi t_k \zeta$. 
\end{proof}

We will now consider the maps $\Hom_{R}(t_k,D)$. The $t_k$'s are called \emph{Eisenbud Operators} and the $\Hom_{R}(t_k,D)$'s are called \emph{Cohomology Operators}. The cohomology operators define an action on the set $\bigoplus\limits_{i\in \mathbb{Z}} \Hom_{R}( C,\Sigma^i D)$ via composition. We can also define operators  $\Hom_{R}(C,s_k)$.  We would like to show that the action of $\Hom_{R}(t_k,D)$ agrees with the action of $\Hom_{R}(C,s_k)$.

\begin{prop2.5}
Consider the graded $R$-module maps $w:\bigoplus\limits_{i\in \mathbb{Z}} \Hom_{R}(\Sigma^{-i}C, D)\rightarrow \bigoplus\limits_{i\in \mathbb{Z}} \Hom_{R}(\Sigma^{-i}C, D),$ $v:\bigoplus\limits_{i\in \mathbb{Z}} \Hom_{R}(C,\Sigma^i D)\rightarrow \bigoplus\limits_{i\in \mathbb{Z}} \Hom_{R}(C,\Sigma^i D),$ and \\
$\phi:\bigoplus\limits_{i\in\mathbb{Z}}\Hom_{R}(\Sigma^{-i} C,D) \rightarrow \bigoplus\limits_{i\in\mathbb{Z}}\Hom_{R}(C,\Sigma^i D)$ defined by 
$$w_j=\Hom_R((\sigma_C^{-j+2})^{-1},D) \circ \Hom_R(t_k,D) \circ \Hom_R(\sigma_{\Sigma^2C}^{-j-2},D),$$ 
$$v_j=\Hom_R(C,\sigma_{\Sigma^2D}^{j-2}) \circ \Hom_R(C,s_k) \circ \Hom_R(C,(\sigma_{D}^{j})^{-1}),$$
and
$$\phi_j(\alpha) = \sigma_{D}^{j} \circ \alpha \circ \sigma_{C}^{-j} \;\;
\text{for any} \;\; \alpha \in \Hom_{R}(\Sigma^{-j} C,D).$$ 
The following diagram commutes up to homotopy.
$$\xymatrix{
	&\bigoplus\limits_{i\in \mathbb{Z}} \Hom_{R}(\Sigma^{-i}C, D) \ar[dl]_{w} \ar[dr]^{\phi}  \\
	\bigoplus\limits_{i\in \mathbb{Z}} \Hom_{R}(\Sigma^{-i}C, D)\ar[dr]^{\phi} & &  \bigoplus\limits_{i\in \mathbb{Z}} \Hom_{R}(C,\Sigma^i D) \ar[dl]_{v} \\
	& \bigoplus\limits_{i\in \mathbb{Z}} \Hom_{R}( C,\Sigma^i D) } $$
\end{prop2.5} 
\begin{proof}
Let $\alpha\in \Hom_{R}(\Sigma^{-j}C, D)$. Proposition 2.2 implies that 
$$\sigma^{j-2}_{\Sigma^2D}\circ s_k \circ (\sigma^{j}_{D})^{-1} \circ \phi(\alpha) \sim \phi(\alpha) \circ \sigma^{-j-2}_{\Sigma^2C} \circ t_k \circ (\sigma^{-j+2}_{C})^{-1} .$$ 
\end{proof}	
Now consider the diagram

$$\xymatrix{& R[\chi_1,\dots,\chi_c] \ar[dl]_{\eta} \ar[dr]^{\nu} & \\
\bigoplus\limits_{i\in \mathbb{Z}}\Hom_R(C,\Sigma^iC) \ar@/_1.5pc/[dr] & & \bigoplus\limits_{i\in \mathbb{Z}}\Hom_R(D,\Sigma^iD) \ar@/^1.5pc/[dl] \\
& \bigoplus\limits_{i\in \mathbb{Z}}\Hom_R(C,\Sigma^iD) &
}$$ 
where $\eta(\chi_k)= t_k$ and $\nu(\chi_k)= s_k$. We can regard the set $\bigoplus\limits_{i\in \mathbb{Z}} \Hom_{R}(C,\Sigma^i D)$ as a left $\bigoplus\limits_{i\in \mathbb{Z}} \Hom_{R}(D,\Sigma^i D)$ right $\bigoplus\limits_{i\in \mathbb{Z}} \Hom_{R}(C,\Sigma^i C)$ graded bimodule with the actions defined by composition. Then $\bigoplus\limits_{i\in \mathbb{Z}} \Hom_{R}(C,\Sigma^i D)$ is an $R[\chi_1,\dots,\chi_c]$-module with the action $\chi_k \cdot \alpha = w(\alpha)$ for any $\alpha\in \bigoplus\limits_{i\in \mathbb{Z}} \Hom_{R}(C,\Sigma^i D)$. The grading is given by $|\chi_k|=2$ for all $k$ and $|r|=0$ for all $r\in R$. We refer to $R[\chi_1,\dots,\chi_c]$ as the \emph{ring of cohomology operators}.  
 
Corollary 2.5 implies that 
 $$\chi_k \cdot \alpha = w(\alpha) \sim \phi\circ v \circ \phi^{-1}(\alpha)$$ 
 and so the action of $R[\chi_1,\dots,\chi_c]$ on $\bigoplus\limits_{i\in \mathbb{Z}} \Hom_{R}(C,\Sigma^i D)$ is independent of the choice of using $C$ or $D$ for the cohomology operators. We would furthermore like this construction to be independent of the choice of basis of $(f_1,\dots,f_c)$. 

\begin{prop2.6}\emph{(Eisenbud \cite[1.7]{Eisenbud:1980})}
Let $Q,Q'$ be local rings and let $f_1,\dots,f_c$ and $f'_1,\dots,f'_{c'}$ be regular sequences in $Q,Q'$ respectively. Consider the complete intersection rings  $R=Q/(f_1,\dots,f_c)$ and $R' = Q'/(f'_1,\dots,f'_{c'})$. Let $\alpha:Q\rightarrow Q'$ be a ring homomorphism such that $\alpha((f_1,\dots,f_c)) \subseteq (f'_1,\dots,f'_{c'})$ with
$$\alpha(f_i)=\sum\limits_{j=1}^{c'}a_{ij}f'_j \;\;\; \text{for some elements} \;\;\; a_{ij}\in Q'.$$
Also let $C\in \Ktac (R)$ with the family of maps $t_i:C\rightarrow \Sigma^2 C$ as defined above. Consider the $R'$-complex $R' \otimes_R C$ with the family of maps $t'_j$. Then $t'_j$ is homotopic to $\sum\limits_{i=1}^c a_{ij}(R' \otimes_R t_i)$.
\end{prop2.6}
If we let $Q=Q'$ and let $f_1,\dots,f_c$ and $f'_1,\dots,f'_{c}$ be two different $Q$-regular sequences that generate the same ideal and we let $\alpha$ be the identity map on $R$, then the above proposition implies that $\sum\limits_{i=1}^c a_{ij}t_i \sim t'_j$ for some $a_{ij}\in Q$. This gives us $\chi'_j\cdot \alpha = \Hom_R(t'_j,D)(\alpha) = \sum\limits_{i=1}^c a_{ij}\Hom_R(t_i,D)(\alpha)=\sum\limits_{i=1}^ca_{ij}\chi_i\cdot \alpha$ for any $\alpha \in \bigoplus\limits_{i\in \mathbb{Z}} \Hom_{R}(C,\Sigma^i D)$, i.e. each $\chi'_j$ is a linear combination of the $\chi_i$'s. 

\section{Support and Rank Varieties} 
We will now translate the definitions of support and rank varieties of modules to the category of totally acyclic complexes. 
\subsection*{Support Varieties}   
   
   Avramov and Buchweitz \cite[2.1]{Avramov:2000} define the \emph{support variety} $V(Q,f,M,N)$ of a pair of $R$-modules $M,N$ to be the zero set of the annihilator of the $\mathbb{k}[\chi_1,\dots,\chi_c]$-module $E = \bigoplus\limits_{i\geq0}\Ext_R^i(M,N) \otimes_R \mathbb{k}$, i.e. $$V(Q,f,M,N)=\{(b_1,\dots,b_c)\in \mathbb{k}^c|\phi(b_1,\dots,b_n)=0 \hspace{0.2cm} \text{for all} \hspace{0.2cm} \phi\in \Ann E\}\cup \{0\}$$
   We will now give a similar definition in $\Ktac(R)$ by replacing $\bigoplus\limits_{i\geq0}\Ext_R^i(M,N)$ by the $R[\chi_1,\dots,\chi_c]$-module $\bigoplus\limits_{i\geq0} \Hom_{R}(C,\Sigma^i D)$, a submodule of $\bigoplus\limits_{i\in \mathbb{Z}} \Hom_{R}(C,\Sigma^i D)$ that we discussed in the previous section.

  \begin{def3.1}
  	Let $R=Q/(f_1,\dots,f_c)$ where $Q$ is a regular local ring with residue field $\mathbb{k}$ and $f_1,\dots,f_c$ is a regular sequence. Also let $S=R[\chi_1,\dots,\chi_c]$ be the ring of cohomology operators defined by $f_1,\dots,f_c$. Then the \emph{support variety} $V(Q,f,C,D)$ of a pair of complexes $C,D\in \Ktac(R)$ is the zero set of the annihilator of the $S\otimes_R \mathbb{k}$-module $E=\bigoplus\limits_{i\geq0} \Hom_{R}(C,\Sigma^iD)\otimes_R \mathbb{k}$. That is,
  	$$V(Q,f,C,D) = \{(b_1,\dots,b_c)\in \mathbb{k}^c|\phi(b_1,\dots,b_n)=0 \hspace{0.2cm} \text{for all} \hspace{0.2cm} \phi\in \Ann E\}\cup \{0\}.$$
  \end{def3.1}
   This definition is especially useful for computing examples by simply decomposing the differentials of $C$ to find the cohomology operators.

\subsection*{Rank Varieties}   
   
   Carlson \cite[4.1]{Carlson} originally defined the \emph{rank variety} $W(M)$ of a module $M$ over a group algebra $\mathbb{k}G\cong \mathbb{k}[[x_1,\dots,x_c]]/(x_1^p,\dots,x_c^p)$ of an elementary abelian $p$-group $G$, where $p=char(\mathbb{k})$, to be the set
    $$W(M)=\{\bar{a}\in \mathbb{k}^c |   M\!\!\downarrow_{\mathbb{k}[l_a]}  \text{is not free}  \}\cup \{0\}$$
    where $\bar{a}$ is the image of $a=(a_1,\dots,a_c)$ with $a_i\in \mathbb{k}G$ and $l_a=a_1x_1+\dots+a_cx_c \in \mathbb{k}G$.   
   By a theorem of Avramov \cite[7.5]{Avramov:1989}, the projective dimension of a module $M$ over the ring $Q_a=\mathbb{k}[[x_1,\dots,x_c]]/(l_a^p)$ is finite if and only if $M$ is free over $\mathbb{k}[l_a]$. Thus 
   $$W(M)=\{\bar{a} \in \mathbb{k}^c | \Ext_{Q_a}^i(M,\mathbb{k})\neq 0 \;\; \text{for infinitely many i}\}\cup \{0\}.$$ 
   Avramov and Buchweitz \cite[2.5]{Avramov:2000} generalize this for a pair of modules $M,N$ over a complete intersection ring $R=Q/(f_1,\dots,f_c)$ by setting
   $$W(Q,f,M,N) = \{\bar{a}\in \mathbb{k}^c | \Ext_{Q_a}^i(M,N)\neq 0 \;\; \text{for infinitely many i}\}\cup \{0\}.$$
  \begin{Remark}
 Avramov and Buchweitz refer to both $V(Q,f,M,N)$ and $W(Q,f,M,N)$ as the support variety after proving they are isomorphic. 
  \end{Remark}
   We would like to give a similar definition for the rank variety of a pair of totally acyclic complexes $C,D$ in terms of the vanishing of $\Hom_{Q_a}$ where $Q_a$ is a complete intersection ring intermediate to $Q$ and $R$. However, totally acyclic $R$-complexes $C,D$ are not totally acyclic $Q_a$-complexes. The following pair of adjoint triangle functors by Bergh, Jorgensen, and Moore will allow us to go back and forth between $\Ktac(R)$ and $\Ktac(Q_a)$. 
   
   \begin{def3.2}(Bergh-Jorgensen-Moore \cite[2.1]{Bergh:2015})
   	Let $R=Q/I$ where $Q$ is a Gorenstein ring and $I$ is an ideal with $\proj\dim_Q R<\infty$, and let $C \in \Ktac(R)$. Then define $T:\Ktac(R)\rightarrow \Ktac(Q)$ by letting $TC \in \Ktac(Q)$ be a complete resolution of $\Im \;d_o^C$ over Q. Given a chain map $\alpha:C \rightarrow C'$ in $\Ktac(R)$, we have the map $u:\Im(d_o^C) \rightarrow \Im(d_o^{C'})$
   	induced by the map $\alpha$. Then $T\alpha:TC \rightarrow TC'$ is the homotopy equivalence class of the comparison map $\tilde{u}:TC \rightarrow TC'$.
   	Define $S:\Ktac(Q)\rightarrow \Ktac(R)$ by $SD = D \otimes_Q R$ and $S\alpha = \alpha\otimes_Q R$.
   	
   \end{def3.2}
   
   \begin{thm3.3}\emph{(Bergh-Jorgensen-Moore \cite[3.1]{Bergh:2015})}
   	The triangle functors $S$ and $T$ form an adjoint pair, that is, they satisfy the following property:
   	for all $C\in \Ktac(R)$ and $D\in \Ktac(Q)$ there exists a bijection
   	$$\Hom_{Q}(D,TC) \to \Hom_{R}(SD,C)$$
   	which is natural in each variable.
   \end{thm3.3}

  \begin{def3.4}
  	 Let $R=Q/(f_1,\dots,f_c)$ and let $C,D\in \Ktac(R)$. Then the \emph{Rank Variety} $W(Q,f,C,D)$ of a pair of complexes $C,D\in \Ktac(R)$ is given by
  	 $$W(Q,f,C,D)=\{\bar{a} \in \mathbb{k}^c  | \Hom_{Q_a}(TC,\Sigma^i TD))\neq 0\hspace{0.1cm} \text{for infinitely many} \; i>0 \}\cup \{0\}$$
  	 where $a=(a_1,\dots,a_c)\in Q^c$ denotes some lifting of $\bar{a}=(\bar{a_1},\dots,\bar{a_c})\in \mathbb{k}^c$ and $Q_a=Q/f_aQ$ with $f_a=a_1f_1+\dots+a_cf_c \in Q$.
  \end{def3.4}

 It is clear that $V(Q,f,C,D)$ is an algebraic variety. However, it is not at all obvious from the definition that $W(Q,f,C,D)$ is a closed set. We will now spend the next section building the tools necessary to show that $V(Q,f,C,D)$ is the same as $W(Q,f,C,D)$. 
  
\section{Finite Generation}
Our main goal in this section is to translate the Gulliksen and  Avramov-Gasharov-Peeva finiteness theorem for $\Ext$ to the setting of totally acyclic chain complexes using the triangulated structure of $\Ktac(R)$. The following generalization of a result by Avramov \cite[2.3]{Avramov:1989} gives a nice condition for showing the finiteness of a graded $A[x]$-module given the finiteness of a graded $A$-module. 
\begin{lem4.1}

Let $A$ be a non-negatively graded noetherian ring and let $F$ be a graded $A$-module and $E$ be a graded $A[x]$-module. Regarding $F$ as a graded $A[x]$-module by the map $s:A[x]\rightarrow A$ where $s(x)=0$, graded in total degree $|ax^n|=|a|+n|x|$ for any homogeneous element $a\in A$. Let there exist graded $A[x]$-module homomorphisms $\psi:E \rightarrow F$ and $\phi:F \rightarrow E$ and an exact sequence of the form
$$\begin{CD}
@>\psi>> F^i @>\phi>> E^{i+|\phi|} @>x>> E^{i+|\phi|+|x|} @>\psi>> F^{i+|\phi|+|x|+|\psi|} @>\phi>>
 \end{CD}.$$ \\
Then $F$ is a finitely generated graded $A$-module if and only if $E$ is a finitely generated graded $A[x]$-module.

\end{lem4.1}
\begin{proof}
First assume $F$ is finitely generated as an $A$-module. Since $F$ is noetherian, $\Im(\psi)$ is finitely generated. Let $\psi(e_1),\dots,\psi(e_m)$ be its generators. Now consider $G = Ae_1+\dots+Ae_m$ as an $A[x]$-graded submodule of $E$. Now let $e\in E$. Then $\psi(e)= a_1\psi(e_1)+\dots+a_m\psi(e_m).$ Thus $e-(a_1e_1+\dots+a_ne_n)\in \Ker(\psi)$, so $E\subseteq G+\Ker(\psi)$. The opposite containment is obvious, thus $E=G+\Ker(\psi)$. But $\Ker(\psi)=\Im(x)=xE$. Thus $E=G+xE$. Iterating yields 
	$$E= \sum_{i=0}^{n} x^i G +x^{n+1}E.$$
For any homogeneous element $e\in E$ and for $n\gg 0$, $e\notin x^{n+1}E$ so we get $E=(A[x]G)$. Thus $E$ is a finitely generated $A[x]$-module.

Now conversely assume that $E$ is a finitely generated $A[x]$-module. Consider the short exact sequence of graded $A[x]$-modules, 
$$0\rightarrow E/\Ker(\psi)\rightarrow F \rightarrow \Im(\phi) \rightarrow 0.$$
Since $E$ is noetherian, both $E/\Ker(\psi)$ and $\Im(\phi)$ are finitely generated, so $F$ is finitely generated as an $A[x]$ module. So for any $f\in F$, $f=p_1f_1+\dots+p_nf_n$ where $p_i\in A[x]$. But $x \in \Ann(F)$, so $p_if_i=a_if_i$ where $a_i\in A$ is the constant term of $p_i$. Thus $F$ is finitely generated as an $A$ module. 
\end{proof} 

\begin{Remark}
	Since $A$ is assumed to be non-negatively graded and $|x|>0$, $F^i,E^i=0$ for $i<0$.
\end{Remark}

We now need to construct such a sequence of $\Hom$ modules. To prove Gulliken's original theorem, Avramov \cite[2.1]{Avramov:1989} uses a complicated homotopy construction by Shamash and Eisenbud. However, we can instead use the triangulated structure of $\Ktac(R)$ to create the long exact sequence. First, we will need a few tools for handling the triangulated structure. 

Given a chain complex $C$ over the ring $R=Q/(f_1,...,f_c)$, we cannot expect that a lifting of $C$ to $Q$ will itself be a complex. However, we can find liftings of mapping cones of the $t_i's$ that are complexes. \\

\begin{lem4.2}
Let $Q$ be a local ring, $R = Q/(f)$ where $f$ is a nonzerodivisor. Also let $C \in \Ktac(R) $ with $\tilde{d}^C_{n-1} \circ \tilde{d}^C_n = \widetilde{t}_{n}f$ so that $t=\sigma^2_C(\widetilde{t}\otimes_Q R)$ is a chain map from $C$ to $\Sigma^2 C$. Then there exists a lifting, $\widetilde{\Cone(t)}$, to $Q$ of $\Cone(t)$ such that the lifting is a complex. Furthermore this complex is exact. 
\end{lem4.2}
\begin{Remark}
For simplicity of notation we will refer to $\widetilde{Cone(t)}$ as $C^{\sharp}$.
\end{Remark}
\begin{proof}
	Consider the mapping cone differential:
$$d_n^{\Cone(t)} = \left(\begin{array}{cc} d_n^{\Sigma C} & 0 \\  t_{n-1} (\sigma^{1}_{C})^{-1}_{n-1} & d^{\Sigma^2 C}_{n} \end{array}\right).$$ \\
Then lift the mapping cone to $Q$ with the differential: 
$${d}_n^{C^{\sharp}} = \left(\begin{array}{cc} \widetilde{d}^{\Sigma C}_{n} & -f(\widetilde{\sigma}^1_C)_{n-2}(\widetilde{\sigma}^2_C)^{-1}_{n-2} \\ (\widetilde{\sigma}^2_C)_{n-3}\widetilde{t}_{n-1} (\widetilde{\sigma}^1_C)^{-1}_{n-1} & \widetilde{d}^{\Sigma^2 C}_{n} \end{array}\right)$$ 

where $\widetilde{\sigma}\otimes_Q R = \sigma$. To see that this lifting is a $Q$-complex, we need to check that the composition of two consecutive differentials is zero. 

$${d}_{n-1}^{C^\sharp} \circ {d}_{n}^{C^\sharp} = \left(\begin{array}{cc} a & b \\  c  & d  \end{array}\right)$$
$$\text{where}$$
$$a = \widetilde{d}^{\Sigma C}_{n-1}\widetilde{d}^{\Sigma C}_n - f (\widetilde{\sigma}^1_C)_{n-3} \widetilde{t}_{n-1} (\widetilde{\sigma}^1_C)^{-1}_{n-1}$$
$$ b =  -\widetilde{d}^{\Sigma C}_{n-1} (\widetilde{\sigma}^1_C)_{n-1}(\widetilde{\sigma}^2_C)^{-1}_{n-1}f-f(\widetilde{\sigma}^1_C)_{n-3}(\widetilde{\sigma}^2_C)^{-1}_{n-3}\widetilde{d}^{\Sigma^2 C}_{n}$$
$$ c = (\widetilde{\sigma}^2_C)_{n-4}\widetilde{t}_{n-2} (\widetilde{\sigma}^1_C)^{-1}_{n-2} \widetilde{d}^{\Sigma C}_n + \widetilde{d}^{\Sigma^2 C}_{n-2} (\widetilde{\sigma}^2_C)_{n-3} \widetilde{t}_{n-1}(\widetilde{\sigma}^1_C)^{-1}_{n-1} $$
$$ d = -f (\widetilde{\sigma}^2_C)_{n-4}\widetilde{t}_{n-2}(\widetilde{\sigma}^2_C)^{-1}_{n-2} + \widetilde{d}^{\Sigma^2 C}_{n-1}\widetilde{d}^{\Sigma^2 C}_{n}.$$ \\

Since $\widetilde{d}^{\Sigma C}_{n-1} \circ \widetilde{d}^{\Sigma C}_n =f (\widetilde{\sigma}^1_C)_{n-3} \widetilde{t}_{n-1} (\widetilde{\sigma}^1_C)^{-1}_{n-1}$ and  $\widetilde{d}^{\Sigma^2 C}_{n-1} \circ \widetilde{d}^{\Sigma^2 C}_n = f (\widetilde{\sigma}^2_C)_{n-4}\widetilde{t}_{n-2}(\widetilde{\sigma}^2_C)^{-1}_{n-2}$, $a$ and $d$ are clearly zero. Furthermore, since $\widetilde{d}_{n-1}^{\Sigma C} =  (-1)(\widetilde{\sigma}^1_{C})_{n-3} \widetilde{d}_{n-2}^C (\widetilde{\sigma}^1_{C})_{n-1}^{-1}$ and $\widetilde{d}_n^{\Sigma^2 C} =  (\widetilde{\sigma}_{C}^2)_{n-3} \widetilde{d}_{n-2}^C (\widetilde{\sigma}^{2}_{C})_{n}^{-1}$ we have
$$ b =  -(\widetilde{\sigma}_{C}^1)_{n-3} \widetilde{d}_{n-2}^C (\widetilde{\sigma}^1_C)^{-1}_{n-1} (\widetilde{\sigma}^1_C)_{n-1} (\widetilde{\sigma}^{2}_{C})_{n}^{-1} + (\widetilde{\sigma}^1_{C})_{n-3} (\widetilde{\sigma}^2_C)^{-1}_{n-3} (\widetilde{\sigma}^2_C)_{n-3} \widetilde{d}_{n-2}^C (\widetilde{\sigma}^{2}_{C})_{n}^{-1}$$ 
$$= -(\widetilde{\sigma}_{C}^1)_{n-3} \widetilde{d}_{n-2}^C (\widetilde{\sigma}^{2}_{C})_{n}^{-1} + (\widetilde{\sigma}^1_{C})_{n-3} \widetilde{d}_{n-2}^C (\widetilde{\sigma}^{2}_{C})_{n}^{-1} = 0 $$
If we multiply $c$ by $f$ we get, 
$$ fc = f\left( (\widetilde{\sigma}^2_C)_{n-4}\widetilde{t}_{n-2} (\widetilde{\sigma}^1_C)^{-1}_{n-2} \widetilde{d}^{\Sigma C}_n + \widetilde{d}^{\Sigma^2 C}_{n-1} (\widetilde{\sigma}^2_C)_{n-3} \widetilde{t}_{n-1}(\widetilde{\sigma}^1_C)^{-1}_{n-1}\right) = $$  $$-(\widetilde{\sigma}^2_C)_{n-4}\widetilde{d}^C_{n-2}\widetilde{d}^C_{n-1}\widetilde{d}^C_n (\widetilde{\sigma}^1_C)^{-1}_{n-1} + (\widetilde{\sigma}^2_C)_{n-4}  \widetilde{d}^C_{n-2}\widetilde{d}^C_{n-1}\widetilde{d}^C_n (\widetilde{\sigma}^1_C)_{n-1}^{-1} = 0.$$  
But $f$ is a non-zerodivisor, thus $c$ is also zero. Thus ${d}_{n-1}^{C^\sharp} \circ {d}_{n}^{C^\sharp} = 0$ and hence the lifting of the cone is a complex. 

Now to see that the complex is exact, let $\widetilde{x}\in \Ker d_n^{C^{\sharp}}$. Then $x=\widetilde{x}+f\Cone(t) \in \Ker d_n^{\Cone(t)}$. We know that $\Cone(t)$ is exact, thus there exists $y\in {\Cone(t)_{n+1}}$ such that $d_n^{\Cone(t)}(y)=x$. Choose $\widetilde{y}\in C^{\sharp}_{n+1}$ such that $y=\widetilde{y} + f\Cone(t)_{n+1}$. Therefore $d_{n+1}^{C^{\sharp}}(\widetilde{y}) - \widetilde{x} \in fC^{\sharp}_n$. So there exists $\widetilde{z} \in C^{\sharp}_{n}$ such that $d_{n+1}^{C^{\sharp}}(\widetilde{y}) - \widetilde{x} = f\widetilde{z}$. Applying $d_{n}^{C^{\sharp}}$ to both sides we get 
$$fd_{n}^{C^{\sharp}}(\widetilde{z}) = d_{n}^{C^{\sharp}}(f\widetilde{z}) = d_{n}^{C^{\sharp}}(d_{n+1}^{C^{\sharp}}(\widetilde{y}) - \widetilde{x}) = 0.$$
Since $f$ is a non-zerodivizor, we can conclude that $d_{n}^{C^{\sharp}}(\widetilde{z}) = 0$, so $\widetilde{z}\in \Ker d_{n}^{C^{\sharp}}$. This implies that
$$ \Ker d_n^{C^1} \subseteq \Im d_{n+1}^{C^{\sharp}} + f \Ker d_n^{C^{\sharp}}.$$
We already know that $\Im d_{n+1}^{C^{\sharp}} \subseteq \Ker d_n^{C^{\sharp}}$, thus by Nakayama's lemma, $\Ker d_n^{C^{\sharp}} =\Im d_{n+1}^{C^{\sharp}}$.
\end{proof}

The above lemma gives a lifting when $R=Q/(f)$. However, we would like to be able to deal with the case where $R=Q/(f_1,\dots,f_c)$. We can apply an induction argument to get a complex over $Q$ in this case. In particular, we can take $R_k=Q/(f_k,\dots,f_c)$ so that $R_k=R_{k-1}/(f_k)$ and so $R=R_1=R_{2}/(f_1)$. Now let $C\in \Ktac(R)$ and lift $C$ to $R_2$ with $(\widetilde{d}^C)^2 = \widetilde{u}_1 f_1$. We know that the complex $\Cone(u_1)\in \Ktac(R)$ lifts to a complex $C^{\sharp}\in \Ktac(R_{2})$. Now we can lift this complex to get
$$(\widetilde{d}^{C^\sharp})^2  = \widetilde{u}_2f_2 \;\; \text{where} \;\;   u_2= \sigma^2_{C^{\sharp}}(\widetilde{u_2} \otimes_Q R) \;\; \text{and} \;\; u_2:C^{\sharp} \rightarrow \Sigma^2 C^{\sharp}.$$ 
But now we can consider $\Cone(u_2) \in \Ktac(R_{2})$ and lift it to $C^{\sharp 2} \in \Ktac(R_{3})$. We can continue this process until we reach $C^{\sharp c} \in \Ktac(Q)$ with maps given by
$${d}_n^{C^{\sharp c}} = \left(\begin{array}{cc} \widetilde{d}^{\Sigma C^{\sharp c-1}}_{n} & -f_c(\widetilde{\sigma}^1_{C^{\sharp c-1}})_{n-2}(\widetilde{\sigma}^2_{C^{\sharp c-1}})^{-1}_{n-2} \\ (\widetilde{\sigma}^2_{C^{\sharp c-1}})_{n-3} (\widetilde{u}_c)_{n} (\widetilde{\sigma}^1_{C^{\sharp c-1}})_{n-1}^{-1} & \widetilde{d}^{\Sigma^2 C^{\sharp c-1}}_{n} \end{array}\right).$$ 
The notation for this iterated lifted cone is somewhat unwieldy. The next lemma will allow us to replace the iterated cone by $TC$.   

\begin{lem4.3}
Let $R=Q/(f_1,\dots,f_c)$ and let $C\in \Ktac(R)$. Also let $C^{\sharp c}$ be defined as above. Then $C^{\sharp c} \simeq TC$. 
\end{lem4.3}
 \begin{proof} First consider the case where $R=Q/(f_1)$ and consider the following $Q$-complex 
 $$
 \begin{CD}
 \dots @>d_3^{C^{\sharp}}>> C^{\sharp}_2 @>d_2^{C^{\sharp}}>> C^{\sharp}_1 @> \left( \widetilde{d}^C_1 (\widetilde{\sigma}^1_C)^{-1}_1 \;\;\; f(\widetilde{\sigma}^2_C)^{-1}_{0} \right)>> \widetilde{C_0} @>d_0^C \circ \pi>> \Im d_0^{C} @>>> 0.
 \end{CD}
 $$
 Where $d_n^{C^{\sharp}}$ is as defined above and $\pi$ is the natural surjection onto $C_0$. This complex agrees with $C^{\sharp}$ to the left. The complex is clearly exact to the left of $C^{\sharp}_1$. To show it is exact at $C^{\sharp}_1$, let $(x_1,x_2)\in \Ker(\widetilde{d}^C_1 (\widetilde{\sigma}^1_C)^{-1}_1 \;\;\; f(\widetilde{\sigma}^2_C)^{-1}_{0})$. Then 
 $$\widetilde{d}_1^C (\widetilde{\sigma}^1_C)^{-1}_1(x_1) + f(\widetilde{\sigma}^2_C)^{-1}_{0}(x_2) = 0$$ 
 which implies $d_1^C(\sigma^1_C)^{-1}_1(x_1)=0$. But $C$ is exact so $(\sigma^1_C)^{-1}_1(x_1)\in \Im d_2^C$ and so there exists $y_1\in (\Sigma C)_2$ such that $$d_2^C(\sigma^1_C)^{-1}_1(y_1)-(\sigma^1_C)^{-1}_1(x_1)=0.$$ Now lifting to $Q$ there exists $y_2 \in (\Sigma^2 C)_2$ such that  
 $$\widetilde{d}_2^C(\widetilde{\sigma}^1_C)^{-1}_1(y_1)-(\widetilde{\sigma}^1_C)^{-1}_1(x_1)=f(\widetilde{\sigma}^2_C)^{-1}_{0}(y_2).$$ 
 Rearranging this and applying $(\widetilde{\sigma}^2_C)_2 \widetilde{d}_1^C $ to both sides we get  
 $$f(\widetilde{\sigma}^2_C)_{2}\widetilde{t}_{1} (\widetilde{\sigma}^1_C)^{-1}_{2}(y_1) +f (\widetilde{\sigma}^2_C)_{2} \widetilde{d}_1^C(\widetilde{\sigma}^1_C)^{-1}_1(y_2)=(\widetilde{\sigma}^2_C)_2\widetilde{d}_1^C(\widetilde{\sigma}^1_C)^{-1}_1(x_1).$$
 
 But we know that $\widetilde{d}_1^C(\widetilde{\sigma}^1_C)^{-1}_1(x_1) = -f(\widetilde{\sigma}^2_C)^{-1}_{0}(x_2)$ so we have
  $$(\widetilde{\sigma}^2_C)_{0}\widetilde{t}_{1} (\widetilde{\sigma}^1_C)^{-1}_{1}(y_1) +\widetilde{d_1^C}(\widetilde{\sigma}^1_C)^{-1}_1(y_2)=x_2.$$ 
 Therefore $(x_1,x_2)\in \Im d_2^{C^{\sharp}}$. 
 
 To see that the complex is exact at $\widetilde{C}_0$, let $\widetilde{x} \in \Ker d_0^{C} \circ \pi$. Then for $\pi(\widetilde{x})=\widetilde{x}+f\widetilde{C}_0$, $x \in \Ker d_0^C$. Since $C$ is an exact complex, there exists $y\in C_1$ such that $d_1^C(y)=x$. Now choose $\widetilde{y} \in \widetilde{C}_1$ such that $y=\widetilde{y}+f\widetilde{C}_1$. Then $\widetilde{d}_1^{C}(\widetilde{y})-\widetilde{x}\in f\widetilde{C}_1$. So there exists $a\in \widetilde{C}_1$ such that $\widetilde{d}_1^C(\widetilde{y})-\widetilde{x}=af$. Therefore $$\widetilde{d}_1^{C}(\widetilde{\sigma}^1_C)^{-1}_1((\widetilde{\sigma}^1_C)_1(\widetilde{y}))-f(\widetilde{\sigma}^2_C)^{-1}_{0}((\widetilde{\sigma}^2_C)_{0}(a)) =\widetilde{x}$$ 
 which shows that $\widetilde{x}\in \Im\left(\widetilde{d}^C_1 (\widetilde{\sigma}^1_C)^{-1}_1 \;\;\; f(\widetilde{\sigma}^2_C)^{-1}_{0}\right)$.
 Thus the complex is exact and hence is a complete resolution of $\Im d_0^C$. So by definition, $C^{\sharp} \simeq TC$. 
 
 Notice if we let $R=Q/(f_1,\dots,f_c)$ and consider $T_i:\Ktac(Q/(f_{i},\dots,f_c)) \rightarrow \Ktac(Q/(f_{i+1},\dots,f_{c}))$, then $T = T_1 \circ \dots \circ T_c$ since $\Im d_0^C \cong \Im d_0^{T_iC}$. Thus in general, $C^{\sharp c} \simeq TC$. 
 \end{proof}
   We now have the tools to prove the main result of this section. 
\begin{thm4.4}
Let $R = Q/(f_1,\dots,f_c)$ where $Q$ is a local noetherian ring with $f_1,\dots,f_c$ a regular sequence and let $C,D \in \Ktac(R)$. Then $\bigoplus\limits_{i\geq0} \Hom_{Q}(TC,\Sigma^iTD)$ is finitely generated as an $R$-module if and only if
$\bigoplus\limits_{i\geq0} \Hom_{R}(C,\Sigma^iD)$ is finitely generated as an $R[\chi_1,...,\chi_c]$-module.
\end{thm4.4}
\begin{proof} First consider the case where $R = Q/(f_1)$ where $f_1$ is a single nonzerodivisor.  Consider the exact triangle 
$$C\rightarrow \Sigma^2 C \rightarrow \Cone(t_1) \rightarrow \Sigma C.$$
Since $\Ktac(R)$ is a triangulated category, we can apply the contravariant $\Hom$ functor to the exact triangle to get the long exact sequence
 \begin{align*}
 \dots \rightarrow &\Hom_R(\Sigma^iC,D) \rightarrow \Hom_R(\Sigma^{i-1}\Cone(t_1),D) \rightarrow \\
&\Hom_{R}(\Sigma^{i+1}C,D) \rightarrow \Hom_{R}(\Sigma^{i-1}C,D) \rightarrow \dots 
 \end{align*}
for any $D\in \Ktac(R)$.
By Lemma 4.2 
$$\Sigma^i\Cone(t_1) \cong S\Sigma^i C^{\sharp} \in \Ktac(R).$$ 
Also, by Theorem 3.3 
$$\Hom_{R}(\Sigma^iSC^{\sharp},D) \cong \Hom_{Q}(\Sigma^iC^{\sharp},TD).$$ 
So we now have the long exact sequence
\begin{align*} \dots \rightarrow &\Hom_R(\Sigma^iC,D) \rightarrow \Hom_Q(\Sigma^{i-1}C^{\sharp},TD) \rightarrow \\ &\Hom_R(\Sigma^{i+1}C,D) \rightarrow \Hom_R(\Sigma^{i-1}C, D) \rightarrow \dots.
\end{align*} 
So lastly we can apply Lemma 4.1 and find that 
$\bigoplus\limits_{i\geq0} \Hom_{Q}(\Sigma^{-i} C^{\sharp},TD)$ is finitely generated as an $R$-module if and only if $\bigoplus\limits_{i\geq0} \Hom_{R}(\Sigma^{-i}C,D)$ is finitely generated as an $R[\chi]$-module. 

Now let the statement hold for $R'=Q'/(f_1,\dots,f_{k-1})$. Let $R'=Q$ and $R=Q/(f_k)$. We can now apply the same argument as above to the exact triangle 
$$C^{\sharp k-1} \rightarrow \Sigma^2 C^{\sharp k-1} \rightarrow \Cone(u_k) \rightarrow \Sigma C^{\sharp k-1}$$
where $u_k$ is the Eisenbud operator of $C^{\sharp k-1}$ for $k>1$ and $u_1=t_1$. This yields that $\bigoplus\limits_{i\geq 0} \Hom_{Q'}(\Sigma^{-i} C^{\sharp k},TD)$ is finitely generated as an $R$-module if and only if
$\bigoplus\limits_{i\geq 0} \Hom_{R}(\Sigma^{-i}C, D)$ is finitely generated as an $R[\chi_1,\dots,\chi_k]$-module. This is equivalent to the statement $\bigoplus\limits_{i\geq 0} \Hom_{Q'}(C^{\sharp k},\Sigma^i TD)$ is finitely generated as an $R$-module if and only if $\bigoplus\limits_{i\geq 0} \Hom_{R}(C,\Sigma^i D)$ is finitely generated as an $R[\chi_1,\dots,\chi_k]$-module. We can now use lemma 4.3 to replace $C^{\sharp k}$ by $TC$ and so the theorem holds.  
\end{proof}

  \section{Avrunin-Scott for $\Ktac(R)$}

 We are now ready to show that the support variety and rank variety of a pair of totally acyclic complexes are the same, translating the Avrunin-Scott theorem \cite[1.1]{Avrunin} to the setting of $\Ktac(R)$. The proof follows the same logic that Avramov and Buchweitz \cite[2.5]{Avramov:2000} use to prove a similar equivalence for the support variety of a pair of modules over a complete intersection ring.

 \begin{thm5.1}
 	Let $R=Q/(f_1,\dots,f_c)$ where $(Q,m,\mathbb{k})$ is a regular local ring and $f=f_1,\dots,f_c$ is a regular sequence. Then $V(Q,f,C,D)= W(Q,f,C,D)$. 
 \end{thm5.1}
 \begin{proof}
 Since $0\in V(Q,f,C,D)$, we can assume that $a_i \neq 0$ for some $i$. Since $f$ is a regular sequence, there exists a regular sequence $f'=f_1',\dots,f'_c$ such that $f$ and $f'$ generate the same ideal and $(a_1,\dots,a_n)=(0,\dots,0,1)$. This means that $f_a=f'_c$, so we only have to show that $\mathbb{k}(0,\dots,0,1)\subseteq V(Q,f,C,D)$ if and only if $\Hom_{Q_a}(TC, \Sigma^i TD)\neq 0$ for infinitely many $i$. We will actually prove the contrapositive of the statement. That is, $\Hom_{Q_a}(TC, \Sigma^i TD)=0$ for $i \gg 0$ if and only if $\mathbb{k}(0,\dots,0,1)\nsubseteq V(Q,f,C,D)$. 

 Let $Q'=Q/(f'_c)$ and recall the notation from definition 3.1, $S=R[\chi_1,\dots,\chi_c]$ and $E=\bigoplus\limits_{i\geq 0} \Hom_R(C,\Sigma^iD)\otimes_R \mathbb{k}$. The sequence $f'_1,\dots,f'_{c-1}$ is regular on $Q'$ so we can consider the cohomology ring of $S'=R[\chi_1',\dots,\chi_{c-1}']$ of $Q'$ defined by this sequence. By proposition 2.6, the action of $\chi'_i$ on $E$ agrees with the action of $\chi_i$ and thus by mapping $\chi_i'$ to $\chi_i$ we can consider $S'$ to be a sub ring of $S$.
 	 
 By Theorem 4.4, we know that $\Hom_{Q'}(TC, \Sigma^iTD)=0$ for $i\gg 0$ if and only if the $S'$-module $\bigoplus\limits_{i\geq 0} \Hom_{R}(C,\Sigma^iD)$ is finitely generated. But by Nakayama's lemma, this happens if and only if $E$ is finitely generated over $S'\otimes_R \mathbb{k}$. This is equivalent to  $E'=E/E(\chi_1',\dots,\chi_{c-1}')\cong E/E(\chi_1,\dots,\chi_{c-1})$ being finitely generated over $(S'\otimes_R \mathbb{k})/(\chi_1',\dots,\chi_{c-1}') \cong \mathbb{k}$. 
 
 Now let $\mathcal{R}=S\otimes_R \mathbb{k}$ and $\mathcal{R}'=\mathcal{R}/(\chi_1,\dots,\chi_{c-1})$. Then $\Rank_{\mathbb{k}} E'$ is finite if and only if $\Supp_{\mathcal{R}}(E') = \{(\chi_1,\dots,\chi_c)\}$. Also since E is finitely generated over $\mathcal{R}$ and $E'=E\otimes_R \mathcal{R}$, we have
 $$\Supp_{\mathcal{R}}(E') = \Supp_{\mathcal{R}} (E \otimes_R \mathcal{R}') = \Supp_{\mathcal{R}}(E)\bigcap \Supp_{\mathcal{R}}(\mathcal{R}').$$
 This is true if and only if
 $$\sqrt{\Ann E}\;\;\bigcap\;\sqrt{\Ann \mathcal{R}'}=(\chi_1,\dots,\chi_c).$$
 Since we have $Z(\Ann_{\mathcal{R}}E)=V(Q,f,C,D)$, $Z(\Ann_{\mathcal{R}}\mathcal{R}')=\mathbb{k}(0,\dots,0,1)$, and $Z(\Ann_{\mathcal{R}}(E'))=0$, the nullstellensatz implies the previous line is true if and only if 
 $$V(Q,f,C,D)\bigcap \mathbb{k}(0,\dots,0,1)=\{0\} \in \mathbb{k}^c.$$
 So looking at the above list of if and only if statements, we now have \\ $\Hom_{Q'}(TC, \Sigma^i TD)=0$ for $i\gg 0$ if and only if $\mathbb{k}(0,\dots,0,1)\nsubseteq V(Q,f,C,D)$.
 \end{proof}

\section{Properties of Support/Rank Varieties}
 We will now prove a few basic results about support/rank varieties in $\Ktac(R)$. The same properties hold for support varieties of modules over a complete intersection ring and the proofs are similar \cite[5.6]{Avramov:2000} cf. \cite[2.2]{Jorgensen:2015}. 
 
 The first theorem provides a nice analogue to Dade's Lemma \cite{Dade} for $\Ktac(R)$. 
  \begin{thm6.1}
     	Let $R=Q/(f_1,\dots,f_c)$ where $(Q,m,\mathbb{k})$ is a regular local ring, $f_1,\dots,f_c$ is a regular sequence, and $C,D \in \Ktac(R)$. Then $\Hom_R(C,\Sigma^iD)=0$ for $i\gg0$ if and only if $V(Q,f,C,D)=\{0\}$ 
    \end{thm6.1}
    \begin{proof}
    Consider the long exact sequence of $\Hom$ from the proof of theorem 4.4. If  $\Hom_{R}(C,\Sigma^nD)$ vanishes for large $n$ then $\Hom_{Q_a}(TC,\Sigma^nTD)$ must also vanish for large $n$. Thus for any $\bar{a} \in \mathbb{k}^c$ such that $\bar{a}\neq0$, $\bar{a}\notin V(Q,f,C,D)$. 
    
    Now conversely let $V(Q,f,C,D)=\{0\}$ and let $E= \bigoplus\limits_{i\geq0} \Hom_{R}(C,\Sigma^iD)\otimes_R \mathbb{k}$. Then $\sqrt{\Ann E}=m$ and so $m^r\subseteq \Ann E$ for some $r>0$. Therefore $m^rE=0$ and so for $n\gg0$ we get $E_{n} \subseteq m^rE_{\leq n-r}=0$. Hence for $n\gg0$, $E_n=\Hom_R(C,\Sigma^nD)=0$. 
    \end{proof}
 
 \begin{thm6.2}
   	Let $R=Q/(f_1,\dots,f_c)$ where $(Q,m,\mathbb{k})$ is a regular local ring, $f_1,\dots,f_c$ is a regular sequence, and $C,D \in \Ktac(R)$. Then the following hold.
    \begin{enumerate}
   \item If $D$ is a complete resolution of the residue field $\mathbb{k}$, then
    $$V(Q,f,D,D)=\mathbb{k}^c.$$
   \item  $V(Q,f,C,D)=V(Q,f,C,C)\cap V(Q,f,D,D)$
   \item If $D$ is a complete resolution of the field $\mathbb{k}$, then $$V(Q,f,C,D)=V(Q,f,C,C).$$
   \item If $C \rightarrow C' \rightarrow \Cone(\alpha) \rightarrow \Sigma C$ and $D\rightarrow D' \rightarrow \Cone(\beta) \rightarrow \Sigma D$ are exact triangles in $\Ktac(R)$, then we have the following inclusions.
   
   \begin{align*}
   V(Q,f,C,D) &\subseteq V(Q,f,C',D)\cup V(Q,f,\Cone(\alpha),D) \\
    V(Q,f,C',D) &\subseteq V(Q,f,C,D) \cup V(Q,f,\Cone(\alpha),D) \\
    V(Q,f,\Cone(\alpha),D) &\subseteq  V(Q,f,C,D) \cup V(Q,f,C',D) \\
    V(Q,f,C,D)&\subseteq V(Q,f,C,D')\cup V(Q,f,C,\Cone(\beta)) \\
    V(Q,f,C,D') &\subseteq V(Q,f,C,D) \cup V(Q,f,C,\Cone(\beta)) \\
    V(Q,f,C,\Cone(\beta)) &\subseteq  V(Q,f,C,D) \cup V(Q,f,C,D') \\
   \end{align*}
 
    \end{enumerate}
    \end{thm6.2}
   
    \begin{proof} \mbox{}
    
    \begin{enumerate}
    
   \item 
   By \cite[2]{Eisenbud:1980}, $\proj\dim_{Q_a} \mathbb{k} = \infty$ thus $\Hom_{Q_a}(TC,\Sigma^n TC)\neq 0$ for infinitely many $n$. Thus $V(Q,f,C,C)=\mathbb{k}^c.$ \\
   
   \item 
    First let $\bar{a}\in V(Q,f,C,C)\cap V(Q,f,D,D)$ but $\bar{a}\notin V(Q,f,C,D)$. Then $\Hom_{Q_a}(TC,\Sigma^n TC)\neq 0$ and $\Hom_{Q_a}(TD,\Sigma^nTD)\neq 0$ for infinitely many $n$ but $\Hom_{Q_a}(TC,\Sigma^nTD)=0$ for $n\gg0$. Then $TC\simeq0\simeq TD$ but $\Hom_{Q_a}(TC,\Sigma^nTD)=0$ for $n\gg0$, a contradiction. Thus $V(Q,f,C,C)\cap V(Q,f,D,D) \subseteq V(Q,f,C,D)$. 
    
    Now let $p(\chi)\in \Ann\bigoplus\limits_{i\geq 0} \Hom_{R}(C,\Sigma^iC)$. Then $p(\Hom_R(t,D))(\alpha)\sim0$ for all $\alpha \in \bigoplus\limits_{i\geq 0} \Hom_{R}(C,\Sigma^iC)$. Hence $p(t)\sim 0$ and so $p(\Hom_R(t,D))(\beta)\sim 0$ for all $\beta \in \bigoplus\limits_{i\geq 0} \Hom_{R}(C,\Sigma^iD)$ and so $p(\chi)\in\Ann\bigoplus\limits_{i\geq 0} \Hom_{R}(C,\Sigma^iD)$. Similarly if $p(\chi)\in \Ann\bigoplus\limits_{i\geq 0} \Hom_{R}(D,\Sigma^iD)$ then $p(\chi)\in\Ann\bigoplus\limits_{i\geq 0} \Hom_{R}(C,\Sigma^iD)$. Thus $$\Ann\bigoplus\limits_{i\geq 0} \Hom_{R}(C,\Sigma^iC) \bigcap \Ann\bigoplus\limits_{i\geq 0} \Hom_{R}(D,\Sigma^iD) \subseteq \Ann\bigoplus\limits_{i\geq 0} \Hom_{R}(C,\Sigma^iD)$$ which implies that $ V(Q,f,C,D) \subseteq V(Q,f,C,C)\cap V(Q,f,D,D)$.\\
        
   \item   
   $$V(Q,f,C,D) = V(Q,f,C,C)\cap V(Q,f,D,D) = $$ 
      $$V(Q,f,C,C)\cap \mathbb{k}^c = V(Q,f,C,C)$$    
      
   \item Applying $\Hom_{Q_a}(\underline{\hspace{.3cm}},D)$ to the exact triangle $$TC \rightarrow TC' \rightarrow T\Cone(\alpha) \rightarrow \Sigma TC$$ we get the long exact sequence
   $$\cdots \rightarrow \Hom_{Q_a}(T\Cone(\alpha),TD)\rightarrow \Hom_{Q_a}(TC',TD) \rightarrow \Hom_{Q_a}(TC,TD)\rightarrow \cdots.$$
   If $\bar{a}\in \mathbb{k}$ is not in $V(Q,f,C',D)\cup V(Q,f,\Cone(\alpha),D)$, then by the rank variety definition $\Hom_{Q_a}(TC',\Sigma^nTD)=\Hom_{Q_a}(T\Cone(\alpha),\Sigma^nTD)=0$ for $n\gg0$. Thus by the long exact sequence, $\Hom_{Q_a}(TC,\Sigma^nTD)=0$ for $n\gg0$ and so $\bar{a}\notin V(Q,f,C,D)$. Hence we have $V(Q,f,C,D)\subseteq V(Q,f,C',D)\cup V(Q,f,\Cone(\alpha),D)$. The other inclusions can be proved similarly. 
     \end{enumerate}
    \end{proof}

\section{Examples}
We will now compute some simple examples of support varieties in $\Ktac(R)$. In the first example we let the field be the complex numbers and choose a simple periodic complex $C$. The resulting support variety $V(Q,f,C,C)$ is a line in $\mathbb{C}^2$.
\begin{example1}

Let $R=\mathbb{C}[[x,y]]/(x^2,y^2)$ and the $C$ be the complex
$$\begin{CD}
C:\cdots @>d^C_{i+1}>> R^2 @>d^C_i>> R^2 @>d^C_{i-1}>>  R^2 @>d^C_{i-2}>> \cdots
\end{CD}$$ 
$$\text{where} \;\;\; d^C_i = \left(\begin{array}{cc} x & y\sqrt{2} \\ y\sqrt{2}  & -x  \end{array}\right).$$

$$\text{Then} \;\;\; \widetilde{d}^2 = \left(\begin{array}{cc} x^2+2y^2 & 0 \\ 0  & x^2+2y^2  \end{array}\right) = x^2\widetilde{t_1}+y^2\widetilde{t_2}$$
 
$$\text{and so} \;\; t_1 = \left(\begin{array}{cc} 1 & 0 \\ 0  & 1  \end{array}\right) \text{and}\;\; t_2 = \left(\begin{array}{cc} 2 & 0 \\ 0  & 2  \end{array}\right)$$
Therefore $$\Ann(\Hom_{R}(C,\Sigma^iC)\otimes_R \mathbb{C}) = (\chi_1-2\chi_2)$$
and so $$V(Q,f,C,C)= Z(\chi_1-2\chi_2) = \{(a,2a)|a\in \mathbb{C} \}.$$				
\end{example1}

The next example illustrates property 2 from theorem 6.2, i.e. that $V(Q,f,C,D)=V(Q,f,C,C)\cap V(Q,f,D,D)$. 
\begin{example2}

Let $R=\mathbb{k}[[x,y,z]]/(x^2,y^2,z^2)$ and let $C$ be the complex 

$$\begin{CD}
C: \cdots @>d^C_{i+1}>> R^2 @>d^C_i>> R^2 @>d^C_{i-1}>>  R^2 @>d^C_{i-2}>> \cdots
\end{CD}$$
\\
where $d^C_i = \left(\begin{array}{cc} x & 0 \\ 0  & y  \end{array}\right)$ for all $i$. 
Then $$(\widetilde{d}^C)^2 = \left(\begin{array}{cc} x^2 & 0 \\ 0  & y^2  \end{array}\right) = x^2\widetilde{t_1}+y^2\widetilde{t_2}+z^2\widetilde{t_3}$$
and so $$t_1 = \left(\begin{array}{cc} 1 & 0 \\ 0  & 0  \end{array}\right), \;t_2 = \left(\begin{array}{cc} 0 & 0 \\ 0  & 1  \end{array}\right), \;\text{and}\;\; t_3 = \left(\begin{array}{cc} 0 & 0 \\ 0  & 0  \end{array}\right).$$
Therefore $$\Ann(\bigoplus\limits_{i\geq 0}\Hom_{R}(C,\Sigma^iC)\otimes_R \mathbb{k}) = (\chi_1\chi_2,\chi_3),$$
thus the support variety of $C$ is $$V(Q,f,C,C) = Z(\chi_1\chi_2,\chi_3)$$ i.e. the union of the $x$ and $y$ axes. 

Now let $D$ be the complex 

$$\begin{CD}
D: \cdots @>d^D_{i+1}>> R^2 @>d^D_i>> R^2 @>d^D_{i-1}>>  R^2 @>d^D_{i-2}>> \cdots
\end{CD}$$
\\
where $d^D_i = \left(\begin{array}{cc} x & 0 \\ 0  & z  \end{array}\right)$ for all $i$. Then $$(\widetilde{d}^D)^2 = \left(\begin{array}{cc} x^2 & 0 \\ 0  & z^2  \end{array}\right) = x^2\widetilde{s_1}+y^2\widetilde{s_2}+z^2\widetilde{s_3}$$
and so $$s_1 = \left(\begin{array}{cc} 1 & 0 \\ 0  & 0  \end{array}\right),\; s_2 = \left(\begin{array}{cc} 0 & 0 \\ 0  & 0  \end{array}\right), \; \text{and}\;\; s_3 = \left(\begin{array}{cc} 0 & 0 \\ 0  & 1  \end{array}\right).$$
Therefore the support variety of $D$ is $$V(Q,f,D,D)=Z(\chi_1\chi_3,\chi_2)$$
i.e. the union of the $x$ and $z$ axes. 

When computing $V(Q,f,C,C)$ we can use $t_k$ and $\Hom_R(t_k,C)$ interchangeably. Similarly we can use $s_k$ and $\Hom_R(s_k,D)$ interchangeably when computing $V(Q,f,D,D)$. However, to compute $V(Q,f,C,D)$ we need to find either $\Hom_R(t_k,D)$ or $\Hom_R(C,s_k)$ for all $1\leq k \leq c$. For any map $h\in \Hom_{R}(C,\Sigma^iD)$, $$h \circ \left(\begin{array}{cc} x & 0 \\ 0  & y  \end{array}\right) = \left(\begin{array}{cc} x & 0 \\ 0  & z  \end{array}\right) \circ h.$$
This implies that $h$ is of the form $h_n=\left(\begin{array}{cc} a_n & 0 \\ 0  & 0  \end{array}\right)$. Since $\Hom_R(t_k,C)(h) = h \circ t_k$ we now have

$$\Hom_{R}(t_1,D)= \left(\begin{array}{cc} 1 & 0 \\ 0  & 0  \end{array}\right) \;\; \text{and} \;\;
 \Hom_{R}(t_2,D)=\left(\begin{array}{cc} 0 & 0 \\ 0  & 0  \end{array}\right)=\Hom_{R}(t_3,D).$$
Therefore $$\Ann(\Hom_{R}(C,\Sigma^iD)\otimes_R \mathbb{k}) = (\chi_2,\chi_3)$$
so we get that $$V(Q,f,C,D) = Z(\chi_2,\chi_3) = V(Q,f,C,C)\cap V(Q,f,D,D)$$ i.e. the support variety of $C$ and $D$ is the $x$ axis. \
\end{example2}

\section*{Acknowledgments}
I would like to thank my advisor, Dr. David Jorgensen, for his help and encouragement in preparing this article. I also thank the algebra seminar group at UTA for the many conversations that helped point me in the right direction.

\bibliographystyle{alpha} 

\begin{thebibliography}{12}
 	\bibitem[1]{Avramov:1989} L. \ Avramov \ \emph{Modules of finite virtual projective dimension}, 1989, Invent. Math., V. 96, pg. 71-101, Springer-Verlag
	\bibitem[2]{Avramov:2000} L. \ Avramov, \ R.O. \ Buchweitz, \emph{Support Varieties and Cohomology over Complete Intersections}, Mar. 2000, Invent. Math., V. 142, pg. 285-315, Springer-Verlag
	\bibitem[3]{Avramov:1997} L. \ Avramov, \ V. \ Gasarov, \ I. \ Peeva, \emph{Complete intersection dimension},1997, Publ. Math. I.H.E.S, 86, 67-114
	\bibitem[4]{Avramov:2002} L. \ Avramov \ A. \ Martsinkovsky  \emph{Absolute, relative, and Tate cohomology of modules of finite Gorenstein dimension}, Proc. London Math. Soc. (3) 85 (2002), 393-440. 
	\bibitem[5]{Avrunin}  G.S. \ Avrunin, \ L. \ Scott,  \emph{Quillen stratification for modules}, 1982, Invent. Math., V.66, pg. 211-286.
	\bibitem[6]{Jorgensen:2015} P. \ Bergh, \ D. \ Jorgensen, \emph{Support Varieties Over Complete Intersections Made Easy}, Sep. 2015, 	arXiv:1509.07828 [math.AC] 
	\bibitem[7]{Bergh:2015} P. \ Bergh, \ D. \ Jorgensen, \ F. \ Moore,  \emph{Totally acyclic approximations}, 2015, In Preparation
	\bibitem[8]{Carlson} J. \ Carlson, \emph{The Varieties and Cohomology Ring of a Module}, Mar. 1982, J. Algebra, V. 85, pg. 104-143
	\bibitem[9]{Dade}E.C.\ Dade, \emph{Endo-permutation modules over $p$-groups. II}, Ann. of Math. (2) V.108, 1978, pg. 317-346
	\bibitem[10]{Eisenbud:1980} D. \ Eisenbud, \emph{Homological Algebra on a Complete Intersection with an Application to Group Representations}, July 1980, Trans. Amer. Math. Soc., V. 260, pg. 35-64
	\bibitem[11]{Gulliksen:1974} T.H. \ Gulliksen, \emph{A Change of Rings Theorem, with applications to Poincare' series and intersection multiplicity}, 1974, Math. Scand. V. 34, pg. 167-183
	\bibitem[12]{Rotman} J. \ J. \ Rotman, \emph{An Introduction to Homological Algebra: Second Edition}, 2009, pg. 340-341, Springer
	
\end{thebibliography}

\end{document}